\documentclass[12pt,draft]{amsart}
\usepackage{amsmath,amsthm,latexsym,amscd,amsbsy,amssymb}
 \relax
\chardef\bslash=`\\ % p. 424, TeXbook

\makeatletter
\def\verbatim{\interlinepenalty\@M \@verbatim
  \leftskip\@totalleftmargin\advance\leftskip2pc
  \frenchspacing\@vobeyspaces \@xverbatim}
\makeatother
\makeatletter
  \def\dgt@k{\dg@DX=-3 \dg@DY=2 \dg@SIZE=3} 
\makeatother

\makeatletter
  \def\dgt@kk{\dg@DX=3 \dg@DY=-1 \dg@SIZE=3}%
\makeatother

\theoremstyle{plain}
\newtheorem{thm}{Theorem}[section]
\newtheorem{cor}[thm]{Corollary}
\newtheorem{lem}[thm]{Lemma}
\newtheorem{pro}[thm]{Proposition}

\theoremstyle{definition}

\newtheorem{defin}[thm]{Definition}

\numberwithin{equation}{section}

\begin{document}

%%%%%%% Begin Topmatter %%%%%%%%%%

\title[Universal $C^{\ast}$-algebra of real rank zero]
{Universal $C^{\ast}$-algebra of real rank zero}
\author{Alex Chigogidze}
\address{Department of Mathematics and Statistics,
University of Saskatche\-wan,
McLean Hall, 106 Wiggins Road, Saskatoon, SK, S7N 5E6,
Canada}
\email{chigogid@math.usask.ca}
\thanks{Author was partially supported by NSERC research grant.}

\keywords{Universal $C^{\ast}$-algebra, real rank, direct limit}
\subjclass{Primary: 46L05; Secondary: 46L85}

%%%%%%% End topmatter %%%%%%%%%

\begin{abstract}{It is well-known that every commutative separable unital
$C^{\ast}$-algebra of real rank zero is a quotient of the $C^{\ast}$-algebra
of all compex continous functions defined on the Cantor cube. We
prove a non-commutative version of this result by showing that the
class of all separable unital $C^{\ast}$-algebras of real rank zero
concides with the class of quotients of a certain separable unital
$C^{\ast}$-algebra of real-rank zero.}
\end{abstract}

\maketitle
\markboth{A.~Chigogidze}{Universal $C^{\ast}$-algebra of
real rank zero}

\section{Introduction}\label{S:intro}
It is a well-known observation (see, for instance,
\cite[Theorem 1.3.15]{eng78}) that the Cantor cube
$\{ 0,1\}^{\omega}$ contains
topological copy of any zero-dimensional metrizable compact space.
By Gelfand's duality and by \cite[Proposition 1.1]{brownped91},
this means that every commutative separable unital
$C^{\ast}$-algebra of real rank zero is a quotient of the $C^{\ast}$-algebra
$C\left( \{ 0,1\}^{\omega}\right)$. Our goal in the present note is
to extend this result to the non-commutative case.

For a given class ${\mathcal C}$ of separable unital
$C^{\ast}$-algebras an
element $Z \in {\mathcal C}$ is said to be universal
(in ${\mathcal C}$)
if every other element $X \in {\mathcal C}$ can be represented as the
image of $Z$ under a surjective unital $\ast$-homomorphism. One
particular method of obtaining universal elements for the class
${\mathcal C}$ is to prove the existence of
${\mathcal C}$-invertible morphism
$p \colon C^{\ast}\left( {\mathbb F}_{\infty}\right)\to Z$.
Here
$C^{\ast}\left( {\mathbb F}_{\infty}\right)$ stands for the group
$C^{\ast}$-algebra
of the free group on countable number of generators and the
${\mathcal C}$-invertibility of $p$ means that for any unital $\ast$-homomorphism
$g \colon C^{\ast}\left( {\mathbb F}_{\infty}\right) \to X$, with
$X \in {\mathcal C}$, there exists a unital $\ast$-homomorphism $h \colon Z \to X$
such that $g = h\circ p$. It is easy to see that in such a situation
$Z$ is an universal element in the class ${\mathcal C}$. Indeed,
every element $X$ of ${\mathcal C}$ (and, in general, every separable unital
$C^{\ast}$-algebra) can be represented as the image of
$C^{\ast}\left( {\mathbb F}_{\infty}\right)$ under a
unital $\ast$-homomorphism
$g \colon C^{\ast}\left( {\mathbb F}_{\infty}\right) \to X$.
By the ${\mathcal C}$-invertibility of $p$ and by the fact
that $X \in {\mathcal C}$, there exists a unital $\ast$-homomorphism
$h \colon Z \to X$ such that $g = h\circ p$. Since $g$ is surjective it follows
from the latter equality that $h$ is surjective as well. 

In the present note we consider the case
${\mathcal C} = {\mathcal R}{\mathcal R}_{0}$, where ${\mathcal R}{\mathcal R}_{0}$
denotes the class of all separable unital $C^{\ast}$-algebras
of real rank zero,
and prove that there indeed exists a ${\mathcal R}{\mathcal R}_{0}$-invertible
morphism
$p \colon C^{\ast}\left( {\mathbb F}_{\infty}\right) \to Z$
such that $Z \in {\mathcal R}{\mathcal R}_{0}$ which, as noted, implies that
$Z$ is universal $C^{\ast}$-algebra in the class ${\mathcal R}{\mathcal R}_{0}$.

%%%%%%%%%%%%%%%%%%%%%%%%%%%%%

\section{Preliminaries}\label{S:pre}
All $C^{\ast}$-algebras below are assumed to be unital and all
$\ast$-homomorphisms between unital $C^{\ast}$-algebras are
also assumed to be
unital. When we refer to a unital $C^{\ast}$-subalgebra of a
unital $C^{\ast}$-algebra we implicitly assume that the
inclusion is a unital $\ast$-homomorphism. The set of all
self-adjoint elements of a $C^{\ast}$-algebra $X$ is denoted
by $X_{sa}$. The density
$d(X)$ of a $C^{\ast}$-algebra $X$ is the minimal cardinality
of dense subspaces (in a purely topological sense) of $X$.
Thus $d(X) \leq \omega$ ($\omega$ denotes the first infinite
cardinal number) means that $X$ is separable. 

%%%%%%%%%%%%%%%%%%%%%%%%%%%%%%%%%%%%%%%%%%%
%%%%%%%%%%%%%%%%%%%%%%%%%%%%%%%%%%%%%%%%%%%%

\subsection{Set-theoretical facts}\label{SS:set}
For the reader's convenience we begin by presenting necessary
set-theoretic facts. Their complete proofs can be found
in \cite{book}.

Let $A$  be a partially ordered {\em directed set} (i.e.
for every two elements  $\alpha ,\beta \in A$  there exists
an element  $\gamma \in A$  such that  $\gamma \geq \alpha$ 
and  $\gamma \geq \beta$). We say that a subset
$A_1 \subseteq A$ of $A$ {\em majorates} another subset
$A_2 \subseteq A$ of $A$ if for each element $\alpha_2 \in A_2$
there exists an element $\alpha_1 \in A_1$ such that
$\alpha_1 \geq \alpha_2$. A subset which majorates $A$
is called {\em cofinal} in $A$. A subset of  $A$  is said to
be a {\em chain} if every two elements of it are comparable.
The symbol $\sup B$ , where  $B \subseteq A$, denotes the
lower upper bound of $B$ (if such an element exists in $A$).
Let now $\tau$ be an infinite cardinal number. A subset $B$
of $A$  is said to be $\tau$-{\em closed} in $A$ if for each chain
$C \subseteq B$, with ${\mid}C{\mid} \leq \tau$, we have
$\sup C \in B$, whenever the element $\sup C$ exists in $A$.
Finally, a directed set $A$ is said to be $\tau$-{\em complete}
if for each chain $C$ of elements of $A$ with
${\mid}C{\mid} \leq \tau$, there exists an element
$\sup C$ in $A$. 

The standard example of a $\tau$-complete set can be obtained
as follows. For an arbitrary set $A$ let $\exp A$ denote, as usual,
the collection of all subsets of $A$. There is a natural partial
order on $\exp A$: $A_1 \geq A_2$ if and only if $A_1 \supseteq A_2$.
With this partial order $\exp A$ becomes a directed set.
If we consider only those subsets of the set $A$ which have
cardinality $\leq \tau$, then the corresponding subcollection
of $\exp A$, denoted by $\exp_{\tau}A$, serves as a basic
example of a $\tau$-complete set.

\begin{pro}\label{P:3.1.1}
Let  $\{ A_{t} : t \in T \}$ be a collection of $\tau$-closed and
cofinal subsets of a $\tau$-complete set $A$. If
$\mid T\mid \leq \tau$, then the intersection
$\cap \{ A_{t}: t \in T \}$ is also cofinal
(in particular, non-empty) and $\tau$-closed in $A$ .
\end{pro}

\begin{cor}\label{C:3.1.2}
For each subset $B$, with  $\mid B \mid \leq \tau$, of a
$\tau$-complete set $A$ there exists an element $\gamma \in A$
such that  $\gamma \geq \beta$  for each  $\beta \in B$ .
\end{cor}

\begin{pro}\label{P:search}
Let  $A$  be a $\tau$-complete set, 
$L \subseteq A^2$, and suppose the following three
conditions are satisfied:
\begin{description}
\item[Existence] For each $\alpha \in A$ there exists
$\beta \in A$  such that  $(\alpha ,\beta ) \in L$.
\item[Majorantness] If  $(\alpha ,\beta ) \in L$  and
$\gamma \geq \beta$, then  $(\alpha ,\gamma ) \in L$.
\item[$\tau$-closeness] Let $\{ \alpha_{t} : t \in T \}$
be a chain in $A$ with $|T| \leq \tau$. If
$(\alpha_{t}, \beta ) \in L$ for some
$\beta \in A$ and each $t \in T$, then
$(\alpha ,\beta ) \in L$ where $\alpha =
\sup \{\alpha_{t} \colon t \in T \}$.
\end{description}
   Then the set $\{ \alpha \in A \colon (\alpha ,\alpha ) \in L\}$
is cofinal and $\tau$-closed in $A$.
\end{pro}

%%%%%%%%%%%%%%%%%%%%%%%%%%%%%%%%%%%%%%%%%%%%
%%%%%%%%%%%%%%%%%%%%%%%%%%%%%%%%%%%%%%%%%%%%

\subsection{Direct $C_{\tau}^{\ast}$-systems of
$C^{\ast}$-algebras}\label{SS:spectral}
Recall (see \cite[Section 1.23]{sakai} for details) that a direct system
${\mathcal S} = \{ X_{\alpha}, i_{\alpha}^{\beta}, A\}$ of
unital $C^{\ast}$-algebras consists of a partially ordered directed
indexing set $A$,
unital $C^{\ast}$-algebras $X_{\alpha}$, $\alpha \in A$, and 
unital $\ast$-homomorphisms
$i_{\alpha}^{\beta} \colon X_{\alpha} \to X_{\beta}$,
defined for each
pair of indexes $\alpha ,\beta \in A$ with $\alpha \leq \beta$,
and satisfying
the condition $i_{\alpha}^{\gamma} =
i_{\beta}^{\gamma}\circ i_{\alpha}^{\beta}$ for
each triple of indexes
$\alpha ,\beta ,\gamma \in A$ with $\alpha \leq \beta \leq \gamma$.
The (inductive) limit of the above direct
system is a unital $C^{\ast}$-algebra which is denoted by
$\varinjlim{\mathcal S}$. For each $\alpha \in A$ there
exists a unital $\ast$-homomorphism
$i_{\alpha} \colon X_{\alpha} \to \varinjlim{\mathcal S}$
which will be called
the $\alpha$-th limit homomorphism of $\mathcal S$. 

If $A^{\prime}$ is a
directed subset of the indexing set $A$, then the subsystem
$\{ X_{\alpha}, i_{\alpha}^{\beta}, A^{\prime}\}$ of
${\mathcal S}$ is denoted ${\mathcal S}|A^{\prime}$.

Below in in Section \ref{S:res} we use the concept of the direct
$C_{\tau}^{\ast}$-system introduced in \cite{chi991}.

\begin{defin}\label{D:smooth}
Let $\tau \geq \omega$
be a cardinal number. A direct system
${\mathcal S} = \{ X_{\alpha}, i_{\alpha}^{\beta}, A\}$ of
unital $C^{\ast}$-algebras $X_{\alpha}$
and unital $\ast$-ho\-mo\-morp\-hisms
$i_{\alpha}^{\beta} \colon X_{\alpha} \to X_{\beta}$ is called a {\em direct
$C_{\tau}^{\ast}$-system} if the
following conditions are satisfied:
\begin{itemize}
\item[(a)]
$A$ is a $\tau$-complete set.
\item[(b)]
Density of $X_{\alpha}$ is at most $\tau$
(i.e. $d(X_{\alpha}) \leq \tau$), $\alpha \in A$.
\item[(c)]
The $\alpha$-th limit homomorphism
$i_{\alpha} \colon X_{\alpha} \to \varinjlim{\mathcal S}$ is
an injective $\ast$-ho\-mo\-mor\-phism for each $\alpha \in A$.
\item[(d)]
If $B = \{ \alpha_{t} \colon t \in T\} $ is a chain of
elements of $A$ with
$|T| \leq \tau$ and $\alpha = \sup B$, then the limit
homomorphism
$\varinjlim\{ i_{\alpha_{t}}^{\alpha} \colon t \in T\}
\colon \varinjlim\left({\mathcal S}|B\right)
\to X_{\alpha}$ is an isomorphism.
\end{itemize}
\end{defin}

\begin{pro}[Proposition 3.2, \cite{chi991}]\label{P:exists}
Let $\tau$ be an infinite cardinal number. Every unital
$C^{\ast}$-algebra $X$
can be represented as the limit of a direct
$C_{\tau}^{\ast}$-system
${\mathcal S}_{X} = \{ X_{\alpha}, i_{\alpha}^{\beta},
A \}$ where the indexing set $A$ coincides with $\exp_{\tau}Y$ for some (any)
dense subset $Y$ of $X$ with $|Y| = d(X)$.
\end{pro}
\begin{proof}
If $d(X) \leq \tau$, then consider the direct $C_{\tau}^{\ast}$-system
${\mathcal S}_{X} = \{ X_{\alpha}, i_{\alpha}^{\beta},
\exp_{\tau}d(X) \}$, where $X_{\alpha} = X$ for each
$\alpha \in \exp_{\tau}d(X)$ and
$i_{\alpha}^{\beta} = \operatorname{id}_{X}$ for each
$\alpha ,\beta \in \exp_{\tau}d(X)$ with
$\alpha \leq \beta$.

If $d(X) > \tau$, then consider any subset $Y$ of $X$ such that
$\operatorname{cl}_{X}Y = X$ and $|Y| = d(X)$. Without loss
of generality we may assume that $Y$ contains the unit of $X$. Each
$\alpha \in \exp_{\tau}d(X)$ can obviously be identified with a subset
(denoted by the same letter $\alpha$) of $Y$ of cardinality $\leq \tau$.
Let $X_{\alpha}$ be the smallest $C^{\ast}$-subalgebra of $X$
containing $\alpha$.
If $\alpha ,\beta \in \exp_{\tau}d(X)$ and $\alpha \leq \beta$, then
$\alpha \subseteq \beta$ (as subsets of $Y$) and consequently
$X_{\alpha} \subseteq X_{\beta}$. This inclusion map is denoted
by $i_{\alpha}^{\beta} \colon X_{\alpha} \to X_{\beta}$. It
is easy to verify
that the collection
${\mathcal S}_{X} = \{ X_{\alpha}, i_{\alpha}^{\beta},
\exp_{\tau}d(X)\}$ is indeed
a direct $C_{\tau}^{\ast}$-system such that
$\varinjlim{\mathcal S}_{X} = X$.
\end{proof}

\begin{lem}[Lemma 3.3, \cite{chi991}]\label{L:strong}
If ${\mathcal S}_{X} = \{ X_{\alpha}, i_{\alpha}^{\beta}, A\}$
is a direct $C_{\tau}^{\ast}$-system, then 
\[ \varinjlim{\mathcal S}_{X} = \bigcup\{
i_{\alpha}(X_{\alpha}) \colon \alpha \in A\} .\]
\end{lem}
\begin{proof}
Clearly $\displaystyle \bigcup\{ i_{\alpha}(X_{\alpha})
\colon \alpha \in A\}$
is dense in $\varinjlim{\mathcal S}_{X}$ (this fact
remains true for arbitrary direct
systems of $C^{\ast}$-algebras). Consequently, for any point
$x \in \varinjlim{\mathcal S}_{X}$ there exists a sequence
$\{ x_{n} \colon n \in \omega\}$, consisting of elements from
$\displaystyle \bigcup\{ i_{\alpha}(X_{\alpha})
\colon \alpha \in A\}$, such that $x = \lim\{ x_{n}
\colon n \in \omega\}$. For each $n \in \omega$ choose an index
$\alpha_{n} \in A$ such that
$x_{n} \in i_{\alpha_{n}}\left( X_{\alpha_{n}}\right)$.
By Corollary \ref{C:3.1.2}, there exists an index $\alpha \in A$
such
that $\alpha \geq \alpha_{n}$ for each $n \in \omega$.
Since $i_{\alpha_{n}} = i_{\alpha}\circ i_{\alpha_{n}}^{\alpha}$,
it follows that
\[ x_{n} \in i_{\alpha_{n}}\left( X_{\alpha_{n}}\right) =
i_{\alpha}\left( i_{\alpha_{n}}^{\alpha}\left(
X_{\alpha_{n}}\right)\right) \subseteq i_{\alpha}\left(
X_{\alpha}\right) \;\;\text{for each}\;\; n \in \omega .\]
Finally, since $i_{\alpha}\left( X_{\alpha}\right)$ is closed
in $\varinjlim{\mathcal S}_{X}$, it follows that 
\[ x = \lim\{ x_{n} \colon n \in \omega\} \in
i_{\alpha}\left( X_{\alpha}\right) .\]
\end{proof}

%%%%%%%%%%%%%%%%%%%%%%%%%%%%%%%%%%%%%%%%%%%%%%%%%%%%%%%%%%%
%%%%%%%%%%%%%%%%%%%%%%%%%%%%%%%%%%%%%%%%%%%%%%%%%%%%%%%%

%%%%%%%%%%%%%%%%%%%%%%%%%%%%%%%%%%%
%%%%%%%%%%%%%%%%%%%%%%%%%%%%%%%%

\section{Existence of an universal separable unital
$C^{\ast}$-algebra of real rank zero}\label{S:res}
The real rank of a unital $C^{\ast}$-algebra $A$, denoted by $RR(A)$,
is defined as follows \cite{brownped91}. We say that $RR(A) \leq n$ if for each
$(n+1)$-tuple $(x_{1},\dots ,x_{n+1})$ of self-adjoint elements
in $A$ and every $\epsilon > 0$,
there exists an $(n+1)$-tuple $(y_{1},\dots ,y_{n+1})$ in $A_{sa}$
such that $\sum y_{k}^{2}$ is invertible and
\[ \left\| \sum_{k=1}^{n+1} (x_{k}-y_{k})^{2}\right\| < \epsilon .\]

Obviously unital $C^{\ast}$-algebras of real rank zero are defined
as those in which every self-adjoint element can be arbitrarily
closely approximated by self-adjoint invertible elements.

\begin{pro}\label{P:rrzero}Let $\tau \geq \omega$ and $X$ be an unital $C^{\ast}$-algebra.
Then the following conditions
are equivalent:
\begin{enumerate}
\item
$RR(X) = 0$.
\item
$X$ can be represented as the direct limit of a direct
$C_{\tau}^{\ast}$-system
$\{ X_{\alpha}, i_{\alpha}^{\beta}, A\}$ satisfying the following properties:
\begin{itemize}
\item[(a)]
The indexing set $A$ is cofinal and $\tau$-closed in the
$\tau$-complete set $\exp_{\tau}Y$
for some (any) dense subset $Y$ of $X$ with $|Y| = d(X)$.
\item[(b)]
$X_{\alpha}$ is a $C^{\ast}$-subalgebra of $X$ such that
$RR(X_{\alpha}) = 0$, $\alpha \in A$.
\end{itemize}
\end{enumerate}
\end{pro}
\begin{proof}
The implication $(2) \Longrightarrow (1)$ follows from
\cite[Proposition 3.1]{brownped91}. 

In order to prove the implication $(1) \Longrightarrow (2)$
we first consider a direct $C_{\tau}^{\ast}$-system
${\mathcal S}_{X} = \{ X_{\alpha}, i_{\alpha}^{\beta}, A\}$
with properties indicated in Proposition \ref{P:exists}. Next
consider the following relation $L \subseteq A^{2}$:

\begin{multline*}
 L = \left\{ (\alpha ,\beta ) \in A^{2} \colon \alpha \leq \beta
\;\text{and for each}\; \epsilon > 0 \; \text{and for each}\;
x \in \left( X_{\alpha}\right)_{sa}\; \text{there}\right.\\
\left.\text{exists an invertible} \; y \in \left( X_{\beta}\right)_{sa}\;
\text{such that}\; \left\| y - x\right\| < \epsilon \right\} .
\end{multline*}

\noindent Let us verify conditions of Proposition \ref{P:search}. 

{\em Existence}. Let $\alpha \in A$ and
$x \in \left( X_{\alpha}\right)_{sa}$. First we prove the following
assertion.

\begin{itemize}
\item[$(\ast )_{\left(\alpha , x , \frac{1}{n}\right)}$]
There exist an index $\beta = \beta\left(\alpha , x,\frac{1}{n}\right)
\in A$, $\beta \geq A$ and an invertible element
$y = y\left(\alpha , x,\frac{1}{n} \right)\in \left( X_{\beta}\right)_{sa}$
such that
$\displaystyle \left\| y - x\right\| < \frac{1}{n}$
\end{itemize}

{\em Proof of $(\ast )_{\left(\alpha , x , \frac{1}{n}\right)}$}.
By (1), $RR(X) = 0$ and consequently there exists
an invertible element $y \in X_{sa}$ such that
$\displaystyle \| y - x\| < \frac{1}{n}$. Since ${\mathcal S}_{X}$
is a direct $C^{\ast}_{\tau}$-system
it follows from Lemma \ref{L:strong} that there exists an index
$\alpha^{\prime} \in A$ such that $y \in X_{\alpha^{\prime}}$.
Since $A$ is a directed set there exists an index $\beta = \beta\left(\alpha , x,\frac{1}{n}\right) \in A$ such that $\alpha ,\alpha^{\prime} \leq \beta$.
This obviously implies that
$X_{\beta} \supseteq X_{\alpha} \cup X_{\alpha^{\prime}}
\supseteq X_{\alpha} \cup \{ y\}$. This finishes the proof of
$(\ast )_{\left(\alpha , x , \frac{1}{n}\right)}$. 

For a given element $x \in \left( X_{\alpha}\right)_{sa}$
consider indeces
$\beta\left(\alpha , x,\frac{1}{n}\right) \in A$, 
satisfying conditions
$(\ast )_{\left(\alpha , x , \frac{1}{n}\right)}$,
$ n = 1,2,\dots$. By Corollary \ref{C:3.1.2},
there exists an element $\beta (\alpha ,x) \in A$ such that
$\beta (\alpha ,x) \geq \beta \left(\alpha, x ,\frac{1}{n}\right)$
for each $n = 1,2,\dots$. Obviously
$X_{\alpha} \subseteq X_{\beta (\alpha ,x)}$, $n = 1,2,\dots$.
Note also that
\begin{itemize}
\item[$(\ast )_{(\alpha ,x)}$]
For any $\epsilon > 0$ there exists an invertible element
$y \in \left( X_{(\alpha ,x)}\right)_{sa}$
such that $\| y -x\| < \epsilon$.
\end{itemize}

Since ${\mathcal S}_{X}$ is a direct $C_{\tau}^{\ast}$-system it
follows that
$d(X_{\alpha}) \leq \tau$. Consequently there exists a subset
$Y_{\alpha} \subseteq \left( X_{\alpha}\right)_{sa}$ such that
$\left( X_{\alpha}\right)_{sa} =
\operatorname{cl}_{\left( X_{\alpha}\right)_{sa}} Y_{\alpha}$ and
$| Y_{\alpha}| = d\left( \left(X_{\alpha}\right)_{sa}\right)
\leq d(X_{\alpha}) \leq \tau$. Let
$Y_{\alpha} = \{ y_{\gamma} \colon \gamma < \tau\}$.
For each $\gamma < \tau$ consider an index
$\beta (\alpha ,y_{\gamma}) \in A$ satisfying condition
$(\ast )_{(\alpha ,y_{\gamma})}$. Since $A$ is a $\tau$-complete set,
we conclude by Corollary \ref{C:3.1.2},
there exists an index $\beta = \beta (\alpha) \in A$ such that
$\beta  \geq \beta (\alpha ,y_{\gamma})$ for each
$\gamma \in \tau$.

We claim that $(\alpha ,\beta ) \in L$. Indeed, let $\epsilon > 0$ and
$x \in \left( X_{\alpha}\right)_{sa}$. Since, by the above construction,
$Y_{\alpha} \subseteq \left( X_{\alpha}\right)_{sa}$ is dense in
$\left( X_{\alpha}\right)_{sa}$, there exists a self-adjoint
element $y_{\gamma} \in Y_{\alpha}$ such that
$\displaystyle \| y_{\gamma} - x\| \leq \frac{\epsilon}{2}$. By condition
$(\ast )_{(\alpha ,y_{\gamma})}$ there exists an invertible element 
$y \in \left( X_{(\alpha ,x)}\right)_{sa}$
such that $\displaystyle \| y -y_{\gamma}\| < \frac{\epsilon}{2}$.
Since
$\beta (\alpha) \geq \beta (\alpha ,y_{\gamma})$, it follows that
$X_{\beta (\alpha ,y_{\gamma})} \subseteq X_{\beta (\alpha )}$.
This guarantees that $y \in \left( X_{\beta (\alpha )}\right)_{sa}$.
It only remains to note that $\| x - y \| < \epsilon$.
Therefore $(\alpha, \beta ) \in L$.

{\em Majorantness}. Let $(\alpha ,\beta ) \in L$,
$\gamma \geq \beta$, $\epsilon > 0$ and
$x \in \left( X_{\alpha}\right)_{sa}$.
Since $(\alpha ,\beta ) \in L$, there exists an invertible
element $y \in \left( X_{\beta}\right)_{sa}$ such that
$\| y-x\| < \epsilon$. Since $\gamma \geq \beta$ it follows
that
$\left( X_{\beta}\right)_{sa} \subseteq \left(
X_{\gamma}\right)_{sa}$
which shows that $y \in \left( X_{\gamma}\right)_{sa}$
and proves that
$(\alpha ,\gamma ) \in L$.

{\em $\tau$-closeness}. Suppose that $\{ \alpha_{t} : t \in T \}$
is a chain of indeces in $A$ with $|T| \leq \tau$. Assume also that
$(\alpha_{t}, \beta ) \in L$ for each $t \in T$
where $\beta \in A$. Our goal is to show that
$(\alpha ,\beta ) \in L$ where $\alpha =
\sup \{\alpha_{t} \colon t \in T \}$. Let $\epsilon > 0$ and
$x \in X_{\alpha}$.
Since ${\mathcal S}_{X}$ is a direct $C^{\ast}_{\tau}$-system
it follows that
$X_{\alpha}$ is the direct limit of the direct system
generated by $C^{\ast}$-subalgebras
$X_{\alpha_{t}}$, $t \in T$, and corresponding inclusion homomorphisms.
Consequently there exist an index $t \in T$ and an element
$x_{t} \in \left( X_{\alpha_{t}}\right)_{sa}$ such that
$\displaystyle \| x - x_{t}\| < \frac{\epsilon}{2}$. Since
$(\alpha_{t},\beta ) \in L$ there exists an invertible element
$y \in \left( X_{\beta}\right)_{sa}$ such that
$\displaystyle \| x_{t} - y\| <\frac{\epsilon}{2}$. Clearly
$\| x-y\| < \epsilon$. This shows that $(\alpha ,\beta ) \in L$.

We are now in position to apply Proposition \ref{P:search} which
guarantees that
the set $A^{\prime} = \{ \alpha \in A \colon (\alpha ,\alpha ) \in L\}$ is cofinal
and $\tau$-closed in $A$.
Note here that 
$(\alpha ,\alpha ) \in L$ precisely when for each $\epsilon > 0$ and for
each element $x \in \left( X_{\alpha}\right)_{sa}$ there exists an
invertible element $y \in \left( X_{\alpha}\right)_{sa}$ such
that $\| x-y\| < \epsilon$. 
This means that the direct $C_{\tau}^{\ast}$-system
$\varinjlim{\mathcal S}_{X}^{\prime} =
\{ X_{\alpha}, i_{\alpha}^{\beta}, A^{\prime}\}$ consists of
$C^{\ast}$-subalgebras of $X$ of real rank zero. Clearly
$\varinjlim{\mathcal S}_{X}^{\prime} = X$.
Proof is completed.
\end{proof}

\begin{cor}\label{C:rrzero}
The following conditions are equivalent for any unital
$C^{\ast}$-algebra $X$:
\begin{enumerate}
\item
$RR(X) = 0$.
\item
$X$ can be represented as the direct limit of a direct
$C_{\omega}^{\ast}$-system
$\{ X_{\alpha}, i_{\alpha}^{\beta}, A\}$ satisfying
the following properties:
\begin{itemize}
\item[(a)]
The indexing set $A$ is cofinal and $\omega$-closed in the
$\omega$-complete set $\exp_{\omega}Y$
for some (any) dense subset $Y$ of $X$ with $|Y| = d(X)$.
\item[(b)]
$X_{\alpha}$ is a separable unital $C^{\ast}$-subalgebra of $X$ such that
$RR(X_{\alpha}) = 0$, $\alpha \in A$.
\end{itemize}
\end{enumerate}
\end{cor}

\begin{lem}\label{L:product}
Let $\{ X_{t} \colon t \in T\}$ be an arbitrary collection
of $C^{\ast}$-algebras and
$RR(X_{t}) = 0$ for each $t \in T$. Then
$\displaystyle RR\left(\prod\{ X_{t} \colon t \in T\} \right) = 0$.
\end{lem}
\begin{proof}
This is an elementary exercise. Let
$\displaystyle x = \{ x_{t} \colon t \in T \} \in
\prod\{ X_{t} \colon t \in T\}$
be a self-adjoint element of the product
$\displaystyle \prod\{ X_{t} \colon t \in T\}$.
This obviously means that $x_{t} = x_{t}^{\ast}$ for each $t \in T$.
Let also $\epsilon > 0$. Since $RR(X_{t}) = 0$, $t \in T$, it
follows that there exists a self-adjoint and invertible element
$y_{t} \in X_{t}$
such that $\displaystyle \| x_{t} - y_{t}\| < \frac{\epsilon}{2}$, $t \in T$.
Obviously $\displaystyle y = \{ y_{t} \colon t \in T \} \in
\prod\{ X_{t} \colon t \in T\}$ and $\| x-y\| < \epsilon$ as required.
\end{proof}

%%%%%%%%%%%%%%%%%%%%%%%%

Next we construct a universal separable unital
$C^{\ast}$-algebra $Z$ of real rank zero. Universal in the
sense that any other
separable unital $C^{\ast}$-algebra with real rank zero is a
quotient of $Z$. We note here that the group $C^{\ast}$-algebra
$C^{\ast}\left( {\mathbb F}_{\infty}\right)$
of the free group on countable number of generators is certainly
universal but its real rank is greater than zero. To see the latter
assume the contrary, i.e.
$RR\left( C^{\ast}\left( {\mathbb F}_{\infty}\right)\right) = 0$. By 
\cite{lin},
$\operatorname{cer}\left( C^{\ast}\left(
{\mathbb F}_{\infty}\right)\right) < 1 +\epsilon$
which contradicts the fact \cite[p.370]{phillips} that
$\operatorname{cer}\left( C^{\ast}\left(
{\mathbb F}_{\infty}\right)\right) = \infty$ (here $\operatorname{cer}$
stands for an exponential rank). Therefore
$RR\left( C^{\ast}\left( {\mathbb F}_{\infty}\right)\right) > 0$.

\begin{thm}\label{T:main}
The class ${\mathcal R}{\mathcal R}_{0}$ of all separable unital
$C^{\ast}$-algebras of real rank zero contains an universal element $Z$.
More formally,
there exists a ${\mathcal R}{\mathcal R}_{0}$-invertible
unital $\ast$-homomorphism
$p \colon C^{\ast}\left( {\mathbb F}_{\infty}\right)
\to Z$ where $Z$ is a separable unital
$C^{\ast}$-algebra such that $RR(Z) = 0$.
\end{thm}
\begin{proof}
Let
${\mathcal A} =
\{ f_{t} \colon C^{\ast}\left( {\mathbb F}_{\infty}\right)
\to X_{t}, t \in T\}$ denote the set
(see \cite[Proposition 2.1]{bla}) of all unital
$\ast$-homomorphisms defined on
$C^{\ast}\left( {\mathbb F}_{\infty}\right)$
such that $RR(X_{t}) = 0$. Next consider the product
$\displaystyle \prod\{ X_{t} \colon t \in T\}$. Since
$RR(X_{t}) = 0$ for each $t \in T$
it follows from Lemma \ref{L:product} that
$\displaystyle RR\left( \prod\{ X_{t} \colon t \in T\}\right) = 0$.
The unital $\ast$-homomorphisms $f_{t}$, $t \in T$, define the unique
unital $\ast$-homomorphism
\[ \displaystyle f \colon C^{\ast}\left( {\mathbb F}_{\infty}\right)
\to \prod\{ X_{t} \colon t \in T\}\]

\noindent such that $\pi_{t} \circ f = f_{t}$ for each $t \in T$ (here
$\displaystyle \pi_{t} \colon \prod\{ X_{t} \colon t \in T\} \to X_{t}$
denotes
the corresponding canonical projection $\ast$-homomorphism).
By Corollary \ref{C:rrzero},
$\displaystyle \prod\{ X_{t} \colon t \in T\}$ can be
represented as the limit of the $C^{\ast}_{\omega}$-system
${\mathcal S} = \{ C_{\alpha}, i_{\alpha}^{\beta}, A\}$ such
that $C_{\alpha}$ is a separable unital $C^{\ast}$-algebra
of real rank zero for each $\alpha \in A$. Suppressing injective unital
$\ast$-homomorphisms $i_{\alpha}^{\beta} \colon C_{\alpha} \to C_{\beta}$ we
for notational simplicity assume that $C_{\alpha}$'s are unital
$C^{\ast}$-subalgebras of $\displaystyle \prod\{ X_{t} \colon t \in T\}$.
Let $\{ a_{n} \colon n \in \omega\}$ be a countable dense subset of
$C^{\ast}\left( {\mathbb F}_{\infty}\right)$. By Lemma \ref{L:strong},
for each $n \in \omega\}$ there exists an index $\alpha_{n} \in A$ such that
$f(a_{n}) \in C_{\alpha_{n}}$. By Corollary \ref{C:3.1.2}, there exists an index
$\alpha \in A$ such that $\alpha \geq \alpha_{n}$ for each $n \in \omega$.
Then 
$f(a_{n}) \in C_{\alpha_{n}} \subseteq C_{\alpha}$ for each $n \in \omega$.
This observation coupled with the continuity of $f$ guarantees that
\[
f\left( C^{\ast}\left( {\mathbb F}_{\infty}\right)\right) =
f\left(\operatorname{cl}\left( \{ a_{n} \colon n \in \omega\right)\right)
\subseteq \operatorname{cl}\left( f\left(\{ a_{n} \colon n \in \omega\}
\right)\right) \subseteq \operatorname{cl}C_{\alpha} = C_{\alpha}.
\]

Let $Z = C_{\alpha}$ and $p$ denote the unital $\ast$-homomorphism $f$
considered as the homomorphism of
$C^{\ast}\left( {\mathbb F}_{\infty}\right)$ into $Z$.
Note that $f = i \circ p$, where
$\displaystyle i \colon Z = C_{\alpha}
\hookrightarrow \prod\{ X_{t} \colon t \in T\}$ stands
for the inclusion.

By construction, $RR(Z) = 0$. Let us show that
$\displaystyle p \colon C^{\ast}\left( {\mathbb F}_{\infty}\right) \to Z$
is ${\mathcal R}{\mathcal R}_{0}$-invertible in the sense of Introduction.
In our situation for any unital $\ast$-homomorphism
$g \colon C^{\ast}\left( {\mathbb F}_{\infty}\right) \to X$, where $X$ is a
separable unital $C^{\ast}$-algebra of real rank zero, we
need to establish the existence of a unital $\ast$-homomorphism
$h \colon Z \to X$ such that $g = h\circ p$. Indeed, by definition
of the set ${\mathcal A}$, we conclude that $g = f_{t}$ for some index $t \in T$
(in particular, $X = X_{t}$ for the same index $t \in T$).
Next observe that $g = f_{t} = \pi_{t} \circ f = \pi_{t} \circ i \circ p$.
This allows us to define the required unital $\ast$-homomorphism
$h \colon Z \to X$ as the composition $h = \pi_{t} \circ i$. This
completes the proof of ${\mathcal R}{\mathcal R}_{0}$-invertibility
of $p$ and as a consequence (see Introduction) of the universality of $Z$.
\end{proof}
%%%%%%%%%%%%%%%%%%%%%%%%%%%%%%%%%

\end{document}